\documentclass[a4paper]{article}
\usepackage{latexsym,amsfonts,amscd,amssymb,enumerate,amstext,amsmath,dsfont}
\usepackage[sc,tiny,center,compact]{titlesec}
\usepackage{amsthm}
\usepackage[all]{xy}
\usepackage{titletoc}
\titlecontents{section}[0pt]{\small\scshape}{\contentslabel{1.5em}}{}{\titlerule*[0.75pc]{.}\contentspage}[\textbullet\ ]

\titlecontents{subsection}[1.5em]{\small\scshape}{\contentslabel{2.3em}}{}{\titlerule*[0.75pc]{.}\contentspage}[\textbullet\ ]
\titlecontents{subsubsection}[3.8em]{\small\scshape}{\contentslabel{3.0em}}{}{\titlerule*[0.75pc]{.}\contentspage}[\textbullet\ ]

\author{M. Grime \\ University of Bristol\\ email: Matt.Grime@bris.ac.uk}
\title{Bousfield Localization in Quotients of Module Categories}
\date{January 2006}
\newtheoremstyle{ordinary}{1ex}{0pt}{}{}{\scshape}{.}{\newline}{}
\theoremstyle{ordinary}
\newtheorem{thm}{Theorem}[section]
\newtheorem{df}[thm]{Definition}
\newtheorem{lem}[thm]{Lemma}
\newtheorem{cor}[thm]{Corollary}

\newtheorem{prop}[thm]{Proposition}
\newtheorem{nonex}[thm]{Nonexample}

\newcommand{\colimit}{\varinjlim}

\begin{document}
\maketitle

\begin{abstract}

We examine various triangulated quotients of the module category of a finite group. We demonstrate that these are not compactly generated by the simple modules and present a modification of Rickard's Idempotent Module construction that accounts for this. When the localizing subcategories are sufficiently nice we give an explicit description of the objects in the Bousfield triangles for modules that are direct limits of sequences of finite dimensional modules in terms of homotopy colimits.

\end{abstract}

\section{Introduction}

Traditionally representation theorists have studied finite dimensional modules, but over the last decade (as exemplified in \cite{idempotent}) infinite dimensional modules have been studied both for their own sake, and for the information they give about finite dimensional modules: \cite{classify}.

One of the principal benefits of allowing infinite dimensional modules is that we now can take arbitrary coproducts of modules, and thus we can bring to bear the machinery of topology such as Bousfield localization and Brown representability. 

We fix $w$ some finite dimensional $kG$-module for a finite group $G$ and consider the class of $w$ projective objects  $\mathfrak{p}(w)$ which are summands of $w\otimes x$ for all $x$ in mod$(kG)$, the infinite dimensional analoge $\mathfrak{P}(w)$ is the class of summands for $w\otimes X$ for all modules $X$. These are thick subcategories, and they constitute the projective objects of exact structures on mod$(kG)$ and Mod$(kG)$. The quotient in either case is a triangulated category, and we denote them by \underline{mod}$_w(kG)$ and \underline{Mod}$_w(kG)$. 

  We start by showing how the relatively stable category fundamentally differs from the usual stable category, mainly by exhibiting a module that is the direct limit of relatively projective objects that is not itself relatively projective\footnote{The example was first suggested to me by Rickard as something where direct limits might not commute with quotients. Sadly the author has not found a general proof to show why this should happen in greater generality, though morally it ought to fail most of the time. 
}. This cannot happen in the ordinary case since projective is equivalent to flat and the direct limits of flat modules are flat.  We are thus  led to a situation where the finite dimensional objects no longer \emph{compactly} generate the relatively stable category \underline{Mod}$_w(kG)$.

  We interpret this failure to generate in terms of Bousfield localization. If we take $\mathcal{L}$ the smallest thick subcategory of \underline{Mod}$_w(kG)$ that contains \underline{mod}$_w(kG)$ then we can localize and prove

{\scshape Theorem}

Let $G$ be a finite group, and $M$  a countably generated module so that $M$ is expressible as the direct limit of a sequence of finitely generated modules. Further, let $w$ be any finite dimensional module and \underline{Mod}$_w(kG)$ the relatively stable category. Then
\[ \mathcal{E}_{\mathcal{L} }(M) \cong_W \text{hocolim}(m_i) \]
for any sequence of modules, $m_i$ with direct limit $M$.

  Until we have a better understanding of what filtered homotopy colimits are over more complex filters then we believe this is the best result this method provides, which of course means we should look for a better method.

\section{The Relatively Stable Category}

We fix the standard representation theoretic materal: $G$ will be some finite group and $k$ a field of characteristic $p$ which divides $|G|$. We will study the module categories mod$(kG)$ and Mod$(kG)$ of finite dimensional and possibly infinite dimensional left $kG$-modules, and we shall adopt the convention that lower case roman letters such as $m$ and $x$  refer to finite dimensional modules, so  something in mod$(kG)$, and upper case to (possibly) infinite dimensional modules (objects in Mod$(kG)$).

Let $w$ be some finite dimensional module, and define $\mathfrak{p}(w)$ to be the $w$-projective modules in mod, and $\mathfrak{P}(w)$ the $w$ projective modules in Mod.  These are the classes of all summands of $w\otimes x$ or $w\otimes X$ for all $x$ in mod or $X$ in Mod respectively. If $x$ is the module Ind$_H^G(k)$ obtained by inducing the trivial module where $H$ is a subgroup of $G$ this agrees with the more traditional notion of projective with respect to a subgroup.

 We wish to study the relatively stable module categories \underline{mod}$_w(kG)$  (or \underline{Mod}$_w(kG)$). The underlying objects are the same as mod$(kG)$ (or Mod$(kG)$), and the morphisms are those in mod$(kG)$ (or Mod$(kG)$) modulo the relation $\alpha \sim \beta \iff \alpha - \beta$ factors through some object in $\mathfrak{p}(w)$ (or $\mathfrak{P}(w)$).
 
 These categories are triangulated: the distinguished triangles come from short exact sequences in the relevant module category that split upon tensoring (over $k$) with $w$ (we will refer to these as $w$-split). See \cite{virtual} or \cite{thesis} for a fuller account of this. We will satisfy ourselves with explaining what the shift functor is, and how to embed a morphism $\alpha:x \to y$ into a triangle in \underline{mod}$_w(kG)$ since this summarizes the essence of the proof. Rather than repeatedly saying ``(or Mod)'' we will just do it in the finite dimensional case; it all passes through without alteration in the infinite dimensional case.
 
 \begin{df} The shift functor (to the left) is written $\Omega_w$ and it is the kernel of the epi
 
 \[ \epsilon :w\otimes w^* \otimes x \to x \]
 if we think of the functor $w\otimes?$ with adjoint $w^*\otimes?$ , then this is  the counit of the adjunction. Traditionally it is  also called the evaluation map. The shift to the right is $\Omega^{-1}_w$ and is the cokernel of
 \[ x \to w\otimes w^* \otimes x \]
 the unit of the adjunction.
 With the functorial view point, trivially we see that these are $w$ split morphisms (in the sense that upon tensoring with $w$ they are split epis and monos). We will call these the \emph{canonical} $w$ projective cover and injective hull.

 \end{df}
 \begin{lem} Suppose that $[\alpha]:y \to z$ is a morphism in \underline{mod}$_w(kG)$, then there is a $w$ split short exact sequence of modules
  \[ 0 \to x \to y' \to z \to 0 \]
 such that $y\cong y'$ in mod$_w(kG)$ and the map from $y'$ to $z$ is equivalent to $[\alpha]$, and thus $[\alpha]$ embeds in a  triangle
 \[ x \to y \to z\]
 which is isomorphic (in the stable category) to the distinguished triangle
 
 \[ x\to y' \to z\]
 \end{lem}

\begin{proof} Lift $[\alpha]$ in the quotient to $\alpha$ in mod$(kG)$ and consider the diagram
\[ \xymatrix{    & y \ar[d]^\alpha \\
		w\otimes w^*\otimes z \ar[r]^{\epsilon}& z}\]
we take the pullback $x$ and form the short exact sequences
\[\xymatrix{  0\ar[r] &{\Omega_w(z)} \ar[r] \ar[d] & x \ar[r] \ar[d]                                                 & y \ar[d]^{\alpha}\ar[r] &0\\
		   0 \ar[r] &{ \Omega_w(z)} \ar[r]          & w\otimes w^* \otimes  z \ar[r]^{\epsilon} & z  \ar[r]&0} \]
Then the $w$ split short exact sequence we need  is
\[0 \to x \to y \oplus (w\otimes w^* \otimes y) \to z \to 0\]
and clearly the middle term is isomorphic to $y$ in the quotient category.

\end{proof}

When $w=kG$ we recover the usual stable category. It was for this category that Rickard translated the language of Bousfield localization to produce the idempotent modules \cite{idempotent}.  We must make some modifications to the construction since the structure of Mod$_w$ is somewhat more mysterious than $\underline{\text{Mod}}(kG)$. For example we know the following about the usual stable category (for any proofs we omit see \cite{idempotent})

\begin{enumerate}

\item Any object in Mod$(kG)$ has a well defined  direct sum decomposition
\[ X \cong (X)_{proj} \oplus(X)_{pf} \]
 into a projective component and a component with no projective summands
 
 \item A morphism $[\alpha]: X \to Y$  between objects with no projective summands is an isomorphism in $\underline{\text{Mod}}(kG)$ iff $\alpha$ is  an isomorphism in Mod$(kG)$
 
\end{enumerate}

For a discussion of the Krull-Schmidt type properties of the relatively stable category we refer the reader to \cite{direct}.   We shall give the proof for the second item, since its failure to generalize in some way contains the essential difference between the two cases.

\begin{proof}
 Suppose  $\alpha$ is an isomorphism in the stable category, and let $s$ be a simple module in the bottom of ker$(\alpha)$. Since the composition of the inclusion of $s$ followed by $\alpha$ is zero in the module category there is a map of triangles in the stable category
 
 \[\xymatrix{  0 \ar[r] & X \ar[r]^{\alpha} & Y \\
 		    s \ar[r]^1 \ar[u]& s\ar[r] \ar[u] & 0 \ar[u]}\]
thus the inclusion factors through the injective hull of $s$, which implies that $X$ has an injective (and hence projective) summand , unless $s=0$ and thus $\alpha$ is monic. Similarly $\alpha$ is epic. The reverse implication is trivial
\end{proof}

We opted not to define \emph{the} injective hull in the relative case though at least for finite dimensional modules it is easy to see there is a well defined hull. We can see this lemma fails (at one of the two possible points it might fail) in even the simplest case when $G= C_2\! \times\! C_2$, $H=C_2$, and $w=\text{Ind}_H^G(k)$ and we do not even need to invoke any infinite dimensional arguments. 

\begin{nonex} With $G = <g,h: g^2=h^2=e, hg=gh>$ and $H=<h>$, the two  non-zero indecomposable  $H$ projective modules are  the free module $kG$ and the two dimensional module $v$ with $G$ action given by

\[ \xymatrix{ k\ar[d]^{1+g}\\ k }\]
and $h$ acts as the identity. All other $H$ projective modules are just direct sums of these objects, since $H$ has finite representation type and induction is an exact functor. The injective hull of the trivial module is  the inclusion into the bottom of $v$. There is a non-$H$-projective  $3$ dimensional module $v'$ again described by its $G$ action it looks like

\[\xymatrix{     k  \ar[dr]_{1+g} &  & k \ar[dl]^{1+h} \\
                                                    &k& }\]
($1+h$ annihilates  the top left and bottom simple module, and $1+g$ annihilates top right and bottom.)             The inclusion of $k$ into $v'$ factors through $v$ but $v$ is not a summand of $v'$. 
\end{nonex}
However, at least in the finite dimensional case, we can prove the relativized version of item 2.

\begin{lem} A morphism  $\alpha:x \to y$ between modules with no relatively projective summands  is invertible in mod$(kG)$ iff it is invertible in \underline{mod}$_w(kG)$
\end{lem}

\begin{proof} Again, one direction remains trivial. For the other direction, suppose that $[\alpha]$ is an isomorphism in the quotient. Let $[\beta]$ be a left inverse in the quotient. By Fitting's Lemma $X$ decomposes as a direct sum 

\[ X = X'\oplus X''\]
with $\beta \alpha$ and isomorphism on $X'$ and nilpotent on $X''$. Since $[\beta \alpha]$ is an isomorphism then $X''$ must be zero  in the quotient contradicting the assumption that $X$ has no projective summands, hence $\alpha$ has a left inverse in the module category. Similarly it must have a right inverse in the module category.
\end{proof}

\section{Bousfield Localizations}

There are two variants of what we now call Bousfield localization for a localizing subcategory $\mathcal{L}$ in a triangulated category $\mathcal{T}$. One puts hypotheses on $\mathcal{T}$ and lets $\mathcal{L}$ by arbitrary the other, the one we shall make use of, places the hypotheses on $\mathcal{L}$.  In \cite{idempotent}, owing to the structure of $\underline{\text{Mod}}(kG)$, the former is used. The hypotheses we allude to refer to compact generation.

\begin{df} An object $c$ in $\mathcal{T}$ is compact if the functor $(c,?)_{\mathcal{T}}$ preserves arbitrary direct sums.
\end{df}

\begin{df} $\mathcal{T}$ is compactly generated if there is a set of compact objects $C$ such that $(c,t)_{\mathcal{T}}=0$ for all $c \in C$ implies $t=0$.
\end{df}

The two variants of Bousfield localization are then
 
\begin{thm}[Bousfield Localization, I]
Let $\mathcal{T}$ be a compactly generated triangulated category with arbitrary direct sums and $\mathcal{L}$ any localizing subcategory. If $t$ is any object in $\mathcal{T}$, then there is a triangle in $\mathcal{T}$

\[ \mathcal{E}_{\mathcal{L}}(t) \to t \to \mathcal{F}_{\mathcal{L}}(t) \]
satisfying 

\begin{itemize}
\item $\mathcal{E}_{\mathcal{L}}(t) \in \mathcal{L}$ 
\item Any morphism from an object in $\mathcal{L}$ factors through
\[ \mathcal{E}_{\mathcal{L}}(t) \to t \]
\item  $(l,\mathcal{F}_{\mathcal{L}}(t)_{\mathcal{T}})=0$ for all $l \in \mathcal{L}$ (it is $\mathcal{L}$ local)
\item The objects  $\mathcal{E}_{\mathcal{L}}(t)$ and  $\mathcal{F}_{\mathcal{L}}(t)$ are universal with respect to these properties, ie any other triangle with terms satisfying these criteria is isomorphic in $\mathcal{T}$

\end{itemize}
\end{thm}
This is the variant we use in the case of $\underline{\text{Mod}}(kG)$: 

 \begin{thm} An object $X$ with no projective summands in $\underline{\text{Mod}}(kG)$ is compact iff it is finite dimensional, further the finite dimensional objects are a compact generating set.
 \end{thm}
\begin{proof} Certainly if a module is finite dimensional it is compact, conversely, if a module $X$ is infinite dimensional at least one simple module $s$ must occur with multiplicity $\Lambda$ for some infinite cardinal  in its top. Then the projection  

\[ X \to \coprod_{\Lambda} s \]
factors through a finite subsum, hence for the other copies of $s$ the projection factors through a projective object, contradicting the assumption $X$ had no projective summands. 

It remains to show the finite dimensional objects generate in this sense. Suppose that $X$ is $C$ local where $C$ is a complete  set of representatives of isomorphism classes of finite dimensional modules. Then the inclusion of any simple module into the bottom factors through a projective, hence $X$ is a direct sum of projectives and is zero in $\underline{\text{Mod}}(kG)$.
\end{proof}

As we  have already seen  proofs of this type have little chance of remaining true in the relative case, and as we shall see in the next section we explicitly construct counterexample to show that it fails in a very elementary case.

Instead there is a variant often referred to as finite Bousfield localization:

\begin{thm}[Bousfield Localization, II]

Let $\mathcal{T}$ be any triangulated category with arbitrary direct sums. If $\mathcal{L}$ is a compactly generated localizing subcategory then the triangle

\[ \mathcal{E}_{\mathcal{L}}(t) \to t \to \mathcal{F}_{\mathcal{L}}(t) \]
with the same properties exists for all $t \in \mathcal{T}$.
\end{thm}

For proofs of either of these statements see \cite{bousfield}

\section{(Homotopy) Colimits}

 In an abelian category with arbitrary coproducts of objects it is well known that we may realize the colimit of a sequence
 
 \[ x_1 \to x_2 \to \ldots \]
 as the cokernel of the monomorphism
 
 \[  (1-s): \oplus x_r  \to \oplus x_r \]
 where we use $s$ for the generic shift map from $x_r$ to $x_{r+1}$.  The triangulated version of this is to use distinguished triangles instead of short exact sequences. 
 \[ \xymatrix{  {\oplus}  x_r \ar[r]^{1-s}& {\oplus} x_r \ar[r] &  {\text{hocolim}} x_r \ar[r] & {\oplus} x_r[1] }\]
The resulting object is called the homotopy colimit of the sequence. The question of when colimits exist for other filters is still open, see eg \cite{neeman} for more details. Note we are using the accepted abuse of language to refer to the homotopy colimit of a sequence; normally the completion of a triangle in a triangulated category is non-canonical, though all completions are isomorphic.

In the case of the ordinary stable category, since short exact sequences and distinguished triangles are interchangeable, the colimit of a sequence is stably isomorphic to its homotopy colimit, however this certainly fails for our relatively stable categories. We shall construct an object that is the direct limit of relatively projective objects in the module category but which is not itself relatively projective. We will use the specific case that $G=H\! \times\! C_2$ over a field of characteristic 2, but the proof can be adapted (see \cite{thesis}) to the situation of $G=H\!\times\! C_p$ for any prime $p$. 

We need a technical result to start.

\begin{prop} Let $k$ be a  field of characteristic $2$ and let $G$ be a finite group of even order that is not of finite representation type, then there is an infinite dimensional module that is not pure projective
\end{prop}

For the proof and explanation of the terms (which we do not need here) see \cite{gnacadja}. In particular this means that there is an infinite dimensional $G$ module   $V$ which is the colimit of finite dimensional modules $v_r$ and such that 

\[ 0 \to \oplus v_r \to \oplus v_r \to V \to 0 \]
is not a split short exact sequence, but each of the \emph{finite} truncations

\[ 0 \to \oplus_{r=1}^n  v_r \to \oplus _{r=1}^{n+1} v_r \to v_{n+1} \to 0\]
is a split short exact sequence. We will use a variation on induction to show that we may construct $kG\times C_2$ modules from the truncations   that are $G$-projective and such that the colimit is not.

We first give the $kG\times C_2$ modules. Start with some short exact sequence of $kG$-modules. 

\[ \xymatrix{ 0 \ar[r] & x \ar[r]^{d_1}& y \ar[r]^{d_2} & z \ar[r]^{d_3} &0 }\]
and define a $kG\times C_2$ action on the vector space $x+y + z$. Let a typical element in $G\times C_2$ be written as $(g,e)$ or  $(g,h)$ for $g \in G$ and $h\neq e \in C_2$ and let them act by $g$ and $g(1+d_i)$ respectively.

\begin{lem} The module given by this construction is $G$-projective iff the short exact sequence of $G$ modules

\[ \xymatrix{ 0 \ar[r] & x \ar[r]^{d_1}& y \ar[r]^{d_2} & z \ar[r]^{d_3} &0 }\]
is split
\end{lem}

\begin{proof} 

Consider $x$ as an $kG\times C_2$ module by extending trivially, and consider the $G$-split surjection  $x\!\!\downarrow^{kG\times C_2}_G \uparrow^{kG\times C_p}_G$ which, if we let $x\!\!\downarrow\uparrow$ be given by $(U,e)+(U,h)$ as the standard vector space sum, has the form $(1,1)$. 
If $y$ is $G$-projective, the map $(1,0): y\to x$, must factor through the natural surjection which is $G$-split by the universality of the counit:

\[\xymatrix{   & &  y \ar@{.>}[dll]^{\theta} \ar[d]^{(1,0)} \\
            {x\downarrow\uparrow} \ar[rr]_{(1,1)}& & x }\]
Let $\theta$ be 

\[\begin{pmatrix} \alpha & \beta & \gamma \\
  \delta & \epsilon & \kappa 
\end{pmatrix}\]
where $\alpha \in \text{Hom}_{kG}(x,x)$ $\beta \in \text{Hom}_{kG}(z,x)$ etc. Then $\theta$
 satisfies 

\[ \begin{pmatrix} 1&0&0\end{pmatrix}=\begin{pmatrix}1&1\end{pmatrix}\begin{pmatrix} \alpha & \beta & \gamma \\ \delta & \epsilon & \kappa \end{pmatrix}\]
ie in char 2, $\alpha+\delta =1$, $\beta=\epsilon$, and $\gamma = \kappa$. And theta must also commute with $h$

\[\begin{pmatrix}0&1\\1&0 \end{pmatrix}\begin{pmatrix} \alpha & \beta & \gamma \\ \delta & \epsilon & \kappa \end{pmatrix} = \begin{pmatrix} \alpha & \beta & \gamma \\ \delta & \epsilon & \kappa \end{pmatrix}\begin{pmatrix} 1&0&0 \\d&1&0 \\ 0&d&1 \end{pmatrix}\]

\[\begin{pmatrix} 1+ \alpha& \beta & \gamma\\ \alpha &\beta &\gamma \end{pmatrix}= \begin{pmatrix} \alpha+\beta d & \beta+\gamma d & \gamma \\ 1+\alpha +\beta d & \beta +\gamma d & \gamma \end{pmatrix}\]
whence $\beta d=1$, and the short exact sequence splits.
\end{proof}

We have now shown that each $kG\times C_2$ module afforded from this construction is $G$-projective for each finite truncation, but that given by  the direct limit it is not $G$-projective.

The sensible reader will take on trust the next statement, skip the proof and move on to the corollary. We may perform a similar construction for any characteristic $p > 0$ and $G$ with Sylow$_p$ subgroup of $p$-rank at least 2. The restriction on the Sylow subgroup means that there is an infinite dimensional module that  is not just a direct sum of finite dimensional modules. Suppose we choose such a $V$ as before that is the direct limit of a sequence of finite dimensional modules. And suppose further that we have a short exact sequence of $G$-modules as before

\[ \xymatrix{ 0 \ar[r] & x \ar[r]^{d_1}& y \ar[r]^{d_2} & z \ar[r]^{d_3} &0 }\]
We can construct a $kG\!\times\! C_p$ module by considering the vector space direct sum

\[x_1+x_2+ \cdots x_{p-1} + y + z_1+ \cdots z_{p-1}\]
where $x\cong x_i$ for all $i$ and $z \cong z_j$ for all $j$. We can define a shift on this vector space direct sum.  Label it $s$ defined such that $s$ sends  $x_i$ to $x_{i+1}$ for $1\leq i\leq p-2$ by acting as some isomorphism and acts as $d_1$ on $x_{p-1}$, it acts on $y$  as $d_2$ and by some isomorphism sending$z_i$ to $z_j$. If $h$ generates $H$, we let $h$ act by $1+s$, and thus a generic element $(g,h^r)$ acts via $g(1+s)^r$.  To make our life easier let us pick bases for each copy of $x$ and $y$ so that the shift acts by $1$.

The same observation that this module is $G$ projective iff the corresponding maps factor holds, since it translates into the following linear algebra.Suppose that we break down $\theta$ etc. as we did in the specific case of $p=2$, then

\[\theta = \left( \begin{array}{cccc} \theta_{1,1} & \theta_{1,2} & \cdots & \theta_{1,2p-1} \\ \vdots & \vdots & \vdots &\vdots \\ \theta_{p,1} & \theta_{p,2}& \cdots & \theta_{p,2p-1}   \end{array} \right)\]
On $U\downarrow\uparrow$ $h$ acts by the $p\times p$ matrix

\[\begin{pmatrix} 0 & 0 & \cdots &0      &     1 \\
                  1 & 0 & 0      &\cdots &     0 \\
                  0 & 1 & 0     & \cdots & 0      \\
                  \vdots&\vdots &\ddots &\ddots& \vdots \\
                  0 & \cdots & 0&1&0 \end{pmatrix} \] 
and, if we suppress the subscripts on the $d_i$, it acts on $M$ as

\[\begin{pmatrix} 1      & 0      & \hdotsfor{6}                                                        & 0 \\
                  1      & 1      & 0           &\hdotsfor{5}                                           & 0 \\
                  0      & \ddots      & \ddots         &\ddots           & \hdotsfor{4}                             & 0 \\
                  \vdots &\ddots  & 1      & 1      & 0     &\hdotsfor{3}                 & 0 \\
                  0      & \hdotsfor{1}& 0      &d           &  1                     &0& \hdotsfor{2}  & 0 \\
                  0      & \hdotsfor{2}         &0           &d      &1&0 &\hdotsfor{1}       & 0 \\
                  \vdots &\hdotsfor{3}          &   0      &1&1  &\ddots            &\vdots  \\
                 \vdots& \hdotsfor{4}                   & \ddots&\ddots& \ddots     &0\\
                  0      & \hdotsfor{5}                            &0  &1                            & 1 \end{pmatrix}\] 
We know that $\theta h=h\theta$ and that $(1,\ldots,1)\theta = (1,0,\ldots,0)$ which gives us a system of equations to solve:

\[\sum_{i=1}^p \theta_{i,1} = 1_U \]
and when $j\neq 1$

\[\sum_{i=1}^p \theta_{i,j}= 0 \]
We only need to consider the entries of $\theta$ for which the second index is less than or equal to $p-1$. (We could consider all of the entries, but that isn't necessary since we only need to prove one of the maps in the short exact sequence splits; considering the other entries would show that the second map splits too.) The relations from commuting with $h$ are:

\[\theta_{r,s}=\theta_{r+1.s}+\theta_{r+1,s+1} \]
where the first index is taken mod $p$, and $s+1\leq p-1$.

We also have the relation

\[ \theta_{1,p}d=\theta_{p,p-1}-\theta_{1,p-1} \]
These relations allow us to find $\theta_{1,p}d$ in terms of the $\theta_{i,1}$.

We need some little combinatorial sublemmas though.

\begin{lem}
For all $k$ and for all $s\leq p-1$

\[\theta_{r,s}= \sum_{i=0}^k \binom{k}{i}\theta_{r-k+i,s-k}\]
\end{lem}
\begin{proof} Induct on $k$,
\[ \sum_{i=0}^k (-1)^i\binom{k}{i}\theta_{r-k+i,s-k}\]

\[ = \sum(-1)^i \binom{k}{i}[\theta_{r-k-1+i,s-k-1}-\theta_{r-k+i,s-k-1}]\]
\[=\sum (-1)^i\binom{k}{i}\theta_{r-k-1+i,s-k-1} + \sum(-1)^{i+1}\binom{k}{i}\theta_{r-k+i,s-(k+1)}\]

\[=\sum (-1)^i[\binom{k}{i}+\binom{k}{i+1}]\theta_{r-(k+1)+i,s-(k+1)}\]
as we were required to show.

In the specific case of $\theta_{p,p-1}$ and $\theta_{1,p-1}$ we see that:

\begin{eqnarray*} \theta_{1,p}d &= & \theta_{p,p-1} - \theta_{1,p-1}\\
& = &\sum_{i=1}^{p-2}(-1)^i\binom{p-2}{i}\theta_{p-(p-2)+i,1} - \sum_{i=1}^{p-2} (-1)^i\binom{p-2}{i}\theta_{1-(p-2)+i,1} \\
& =&  \sum_{i=1}^{p-2}(-1)^i\binom{p-2}{i}\theta_{2+i,1} - \sum_{i=1}^{p-2} (-1)^i\binom{p-2}{i}\theta_{3-p+i,1}
\end{eqnarray*}
And since the first index is mod $p$, this is

\[   \sum_{i=1}^{p-2}(-1)^i\binom{p-2}{i}\theta_{2+i,1} - \sum_{i=1}^{p-2} (-1)^i\binom{p-2}{i}\theta_{3+i,1}\]
shifting the index in the second sum, and combining, this is the same as

\[ \sum(-1)^i[\binom{p-2}{i}+\binom{p-2}{i+1}]\theta_{2+1,1}\]
\[ =\sum (-1)^i\binom{p-1}{i}\theta_{2+i.1}\]
Claim: $\binom{p-1}{i} \cong (-1)^i \ \text{mod} p$. And then
\[ \theta_{1,p}d = \sum_i \theta_{i,1} = 1_U \]

and $d$ splits.\end{proof}

Proof of claim:
Recall that $\binom{p}{i}$ is zero mod $p$ for $i$ not equal to $0$ or $p$, and the claim follows from considering the sum of the coefficients $\binom{p-1}{i}$ and $\binom{p-1}{i+1}$.

\begin{cor} The set of finite dimensional modules, whilst compact, does not necessarily  compactly generate \underline{mod}$_w(kG)$ in general
\end{cor}
\begin{proof}
We have shown there is a module $X= \colimit x_r$ where each $x_r$ is relatively projective but $X$ is not.  If $c$ is any finite dimensional module any (lift of) a morphism in the quotient from $c$ to $X$ in the factors through some $x_n$ and is zero in the relatively stable category.
\end{proof}

\section{Localizing with respect to finite dimensional objects}

We suppose that $k$ is a field of characteristic $p$ and that $M$ is a  module that  is the colimit (in the module category) of a sequence of finite dimensional modules, or equivalently that  dim$(M) \leq \aleph_0$. Suppose that 

\[ M \cong \colimit m_i \]
so it fits into a short exact sequence
\[ \xymatrix { 0 \to {\oplus} m_i \ar[r]^{1-s} & {\oplus} m_i \ar[r]^t & M \ar[r] & 0}\]
where we use $s$ to denote the generic inclusion $s_i:m_i \to m_{i+1}$. Since $t(1-s)=0$ in the module category it follows that there is a map of distinguished triangles in the relatively stable category

\[\xymatrix{ {\oplus}m_i \ar[r]^{[1-s]} \ar[d]& {\oplus}m_i \ar[r] \ar@{=}[d]^{[t]} & {\text{hocolim}}m_i \ar[r] \ar[d]^{[\lambda]}& {\oplus (m_i)[1]} \ar[d] \\
  0 \ar[r]& M \ar@{=}[r] & M \ar[r]& 0 }\]
We aim to prove that if we let $\mathcal{S}$ be the smallest thick triangulated subcategory of \underline{Mod}$_w(kG)$ containing \underline{mod}$_w(kG)$ and closed under all direct sums (which is compactly generated) then the first object arising in the Bousfield localization triangle

\[ \mathcal{E}_{\mathcal{L}}(M) \to M\to \mathcal{F}_{\mathcal{L}}(M) \]
is (stably) isomorphic to hocolim$(m_i)$ where $M \cong \colimit(m_i)$, and further that the map  

\[ \mathcal{E}_{\mathcal{L}}(M) \to M \]
may be taken to be $[\lambda]$. Naturally, it is easiest to show that the triangle 

\[ \text{hocolim}(m_i) \to M \to M'\]
obtained by completing the morphism $[\lambda]$ to a triangle in the relatively stable category has the correct universal characterization.  To do this it suffices to show that every morphism in the quotient from a finite dimensional object factors uniquely through the homotopy colimit, for then $M'$ is the universal $\mathcal{L}$ local object and it is indeed a Bousfield triangle. Obviously the only difficulty is showing uniqueness since every map from a finite dimensional object factors through some $m_r$.

First, we need a result about how we may choose lifts of maps from the stable to the module category.

\begin{lem} Suppose that $\beta:x \to M$ is a module map, and that it factors through a relatively projective module, ie $[\beta]=0$. Let $\beta$ factor as $t\beta'$ in the module category, then there is a map $\gamma$ such that $[\beta']=[\gamma]$ and $t\gamma=0$.
\end{lem}

\begin{proof}Since $[\beta]=0$, it factors through $I(x):=w\otimes w^* x$ which is finite dimensional and there is a diagram 
\[\xymatrix{   x \ar[r]^{\iota} &  I(x) \ar[d]\ar[r] & M \\
                                & {\oplus}m_i \ar[ur]^{t} & }\]
Let $\delta$ be the composite $ x \to I(X) \to \oplus m_i$, and set $\gamma=\beta'-\delta$, and we are done.\end{proof}

We are now in a position to prove 

\begin{thm}
The map 

\[\xymatrix{ {\lambda_*}:(x,\text{hocolim}(m_i))_W \ar[r] & (x,M)_W }\]
is an isomorphism.
\end{thm}

\begin{proof} We have observed  it is an epimorphism already. Suppose, now, that $[\alpha]$ is in the kernel of $\lambda_*$, ie $[\lambda\alpha]=0$. In the relatively stable category $[\alpha]$ factors as $[t'\beta]$. Now, $[t\beta]:x \to M$ is zero, in the relatively stable category. By the preceding lemma we may suppose that we have chosen $\beta$ such $t\beta$ is zero in the module category. That is in the module category $\beta=\zeta(1-s)$ for some $\zeta$. 

That is we have

\[ [\alpha] = [t'\beta]=[t'(1-s)\zeta]\]
which is zero since $[t'(1-s)]$ is the composition of two consecutive maps in a distinguished triangle. \end{proof}

We are now in a position to prove the main result of this section.

\begin{thm}

Let $G$ be a finite group, and $k$ an algebraically closed field. Further, let $W$ be any finite dimensional module and \underline{Mod}$_W(kG)$ the relatively stable category, and suppose that $M$ is a module such that dim$_k(M)$ is countable. Then
\[ \mathcal{E}_{\mathcal{L} }(M) \cong_W \text{hocolim}(m_i) \]
where the $m_i$ are finite dimensional submodules whose direct limit is $M$.
\end{thm}
\begin{proof}

In the distinguished triangle

\[ \text{hocolim}(m_i) \to M \to M' \to \text{hocolim}(m_i)[1]\]
every map in the relatively stable category from a finite dimensional module to $M$ factors through the homotopy colimit, and the factoring is unique so that $M'$ is local with respect to the finitely generated modules. These facts allow us to conclude that 

\[ \mathcal{E}_{\mathcal{L}}M \cong_W \text{hocolim}(m_i) \]
\end{proof}

\section{Acknowledgements}

This work was carried out whilst a PhD student at the University of Bristol under the kind and knowledgeable supervision of Jeremy Rickard to whom I am greatly indebted. The position was funded by the EPSRC. I also benefited from the advice and encouragement of Joe Chuang and Thorsten Holm.

\nocite{alperin}
\nocite{benson1}
\nocite{neeman}
\nocite{virtual}
\nocite{ideals}
\nocite{idempotent}
\nocite{okuyama}

\nocite{bousfield}
\addcontentsline{toc}{chapter}{Bibliography}
\bibliographystyle{alpha}

\bibliography{bibliography}

\end{document}